\documentclass[
oneside]{amsart}

\usepackage{inputenc}
\usepackage{cite}

\usepackage{amsthm}

\usepackage{stackengine}[2013-09-11]

 \usepackage{amsmath}
\usepackage{amsfonts}
\usepackage{amssymb}
\usepackage{amsthm}
\usepackage{colonequals}
\usepackage[usenames]{xcolor}

\usepackage[hyperfootnotes=false]{hyperref}

\title[]{A structure theorem for neighborhoods of  compact complex manifolds}
\author[]{Xianghong Gong$^{\dag}$}

\address{Department of Mathematics,
	University of Wisconsin-Madison, Madison, WI 53706, U.S.A.}
\email{gong@math.wisc.edu}

\author{Laurent Stolovitch$^{\dag\dag}$}
\address{CNRS and Laboratoire J.-A. Dieudonn\'e
	U.M.R. 7351, Universit\'e C\^ote d'Azur, Parc Valrose
	06108 Nice Cedex 02, France}
\email{stolo@unice.fr}
\thanks{$^{\dag}$Partially supported by a grant from the Simons
	Foundation (award number: 505027) and NSF grant DMS-2054989. $^{\dag\dag}$Research of L. Stolovitch has been supported by the French government through the UCAJEDI Investments in the Future project managed by the National Research Agency (ANR) with the reference number ANR-15-IDEX-01.  }


\newtheorem{thm}{Theorem}[section]
\newtheorem{cor}[thm]{Corollary}
\newtheorem{prop}[thm]{Proposition}
\newtheorem{lemma}[thm]{Lemma}

\theoremstyle{definition}

\newtheorem{defn}[thm]{Definition}

\newtheorem{rem}[thm]{Remark}

\renewcommand{\th}[1]{\begin{thm}\label{#1}}
	\newcommand{\co}[1]{\begin{cor}\label{#1}}
		\newcommand{\eco}{\end{cor}}
	\renewcommand{\le}[1]{\begin{lemma}\label{#1}}
		\newcommand{\ele}{\end{lemma}}
	\newcommand{\pr}[1]{\begin{prop}\label{#1}}
		\newcommand{\epr}{\end{prop}}

	\newcommand{\ga}{\begin{gather}}
	\newcommand{\ega}{\end{gather}}
	\newcommand{\gan}{\begin{gather*}}
	\newcommand{\egan}{\end{gather*}}
	\newcommand{\al}{\begin{align}}
	\newcommand{\eal}{\end{align}}
	\newcommand{\aln}{\begin{align*}}
	\newcommand{\ealn}{\end{align*}}
	\newcommand{\eq}[1]{\begin{equation}\label{#1}}
	\newcommand{\eeq}{\end{equation}}

	\newcommand{\ci}{~\cite}

\newcommand{\aut}{\operatorname{Aut}}

\newcommand{\C}{{\mathbb C}}

\newcommand{\tcm}[1]{T_CM\otimes S^{#1}N_C^*}

	\newcommand{\codim}{\operatorname{codim}}

	\newcommand{\rank}{\operatorname{rank}}

	\newcommand{\cL}{\mathcal}

\newcommand{\Mt}{\mathcal M}

	\newcommand{\all}{\alpha}

	\newcommand{\del}{\delta}
	\newcommand{\Del}{\Delta}
	\newcommand{\var}{\varphi}
	\newcommand{\e}{\epsilon}
	\newcommand{\om}{\omega}
	\newcommand{\Om}{\Omega}

	\newcommand{\pd}{\partial}

	\newcommand{\re}[1]{(\ref{#1})}
	\newcommand{\rea}[1]{$(\ref{#1})$}
	\newcommand{\rl}[1]{Lemma~\ref{#1}}
	
	\newcommand{\rp}[1]{Proposition~\ref{#1}}
	\newcommand{\rt}[1]{Theorem~\ref{#1}}
	\newcommand{\rd}[1]{Definition~\ref{#1}}
	\newcommand{\rrem}[1]{Remark~\ref{#1}}

	\newcommand{\rrema}[1]{Remark~$\ref{#1}$}

	\newcounter{pp}
	\newcommand{\bpp}{\begin{list}{$\hspace{-1em}(\alph{pp})$}{\usecounter{pp}}}
		\newcommand{\epp}{\end{list}}
	
	\newcounter{ppp}
	\newcommand{\bppp}{\begin{list}{$\hspace{-1em}(\roman{ppp})$}{\usecounter{ppp}}}
		\newcommand{\eppp}{\end{list}}
	
	\def\beq{\begin{equation}}
	\def\eeq{\end{equation}}

\begin{document}

\begin{abstract}We construct an injective map from the set of holomorphic equivalence classes of  neighborhoods $M$ of a compact complex manifold $C$ into $\C^m$ for some $m<\infty$ when  $(TM)|_C$ is fixed and the normal bundle of $C$ in $M$  is either \emph{weakly negative} or $2$-\emph{positive}.
\end{abstract}

\date{\today}
 \maketitle


\setcounter{thm}{0}\setcounter{equation}{0}
	\section{Introduction}
Let $C$ be a  compact complex manifold.   We say that two holomorphic embeddings $f\colon C \hookrightarrow M$ and $\tilde f\colon C  \hookrightarrow\tilde M$ are {\it holomorphic equivalent} if there is a biholomorphic mapping  $F$ from a neighborhood of
$f(C)$ in $M$ into 
a neighborhood of
$\tilde f(C)$ in $\tilde M$ such that $F f=\tilde f$.  To classify   such neighborhoods $M$, we identity $C$ with $f(C)$ via $f$. We also fix  the normal bundle $N_C$ of $C$ in $M$ and $T_CM$, the restriction of $TM$ on $C$. To determine the holomorphic equivalence of two neighborhoods, Grauert~\ci{MR0137127} introduced
the \emph{formal principle}  that
asserts   two neighborhoods of $C$ are  holomorphic equivalent  if they are formally equivalent.	The formal principle holds for  neighborhoods when $N_C$ has nice geometry properties. When $N_C$ is \emph{weakly negative}, Grauert~\ci{MR0137127} proved that the formal principle holds when $C$ has codimension one. Hironaka and Rossi~\cite{MR0171784} extended Grauert's result to   arbitrary positive codimension as well as to the reduced complex spaces $C$ that are exceptional in the sense of Grauert.  Griffiths~\cite{MR206980}  showed that the formal principle holds when $N_C$ is \emph{sufficiently positive} and $\dim C\geq3$. This result was improved by Hirschowitz~\cite{MR621013} including $\dim C\geq2$ with \emph{weakly positive} $N_C$, 
and then by Commichau-Grauert~\ci{MR627752} for $1$-\emph{positive} normal bundle $N_C$ for $\dim C\geq1$.

Despite all these positive results on the formal principle, it remains unknown until now if the equivalence classes of neighborhoods of $C$ with a fixed  $T_CM$ form a \emph{finite} dimensional space
 when the formal principle holds.
The main result of this paper provides a structure for the neighborhoods of $C$ as follows.

\begin{thm}\label{m-finite} Let $C$ be a compact complex manifold.
 Assume that    $N_C$ is either weakly  negative or  $2$-positive.
There is an injective mapping from the set of  holomorphic equivalence classes of neighborhoods of $C$ into  the finite-dimensional space
$$
 \cL H^1(T_CM):=\bigoplus_{\ell\geq2} H^1(C,T_CM\otimes S^\ell N_C^*),
$$
where $S^\ell N_C^*$ is the $\ell$-th symmetric power of the dual bundle $N_C^*$ of $N_C$.
\end{thm}

Note that  $\dim\cL H^1(T_CM)<\infty$ follows from the assumption on $N_C$.
Indeed, by definitions in~\cite[Def.~1, p.~342]{MR0137127}  and \ci[Def.~1, p.~108]{MR627752},  $N_C$ is \emph{weakly negative} (resp. $q$-positive) if the zero section of $N_C$  admits a tubular neighborhood
$Tub(C)$
 in $N_C$ whose boundary is strictly pseudoconvex (resp. of at least
 $q$ negative Levi eigenvalues).
 It is known that $\dim \cL H^1(T_CM)$ is finite, when the zero section of $N_C$ admits a relatively compact neighborhood $W$ in $N_C$ that is either weakly negative or $2$-positive.

The formal principle does not rule out the existence of  an
infinite dimensional moduli space, i.e. equivalence classes.
 In fact, given a (compact and smooth) Riemann surface $C$ and a \emph{positive} line bundle $N_C$, Morrow-Rossi~\ci[p.~323]{MR627765} constructed a
complete set of  the equivalence classes of holomorphic transverse foliations of such $C$, under the smaller group that preserves the foliations,   showing that the moduli space for the equivalent classes is infinite-dimensional. See \rp{realization-MR} for details.

For other closely related results, we mention that when  $C$ is the Riemann sphere and $\rank N_C=1$,   Hurtubise-Kamran~\ci{MR1152943} showed that  the space of  neighborhoods with $1$-positive $N_C$ is infinite-dimensional, and  Mishustin~\ci{MR1250979} constructed a normal form with infinite-dimensional invariants. 
See  Ilyashenko~\ci{MR704627} when $C$ is elliptic curve and $N_C$ is a positive line bundle
 and recent results of
  Falla Luza and Loray~\ci{MR3488114} on the Riemann sphere.
   On the other hand, when $C$ is an elliptic curve embedded in a complex surface with $N_C$ topologically trivial, Arnol'd~\ci{MR0431285} showed that the formal principle
   holds under an extra condition, namely a suitable {\it Diophantine condition} on $T_CM$ when it splits as the direct sum   $TC\oplus N_C$. Doing so, Arnol'd~ showed for the first time that the formal principle \emph{fails} when $N_C$ is topologically trivial for an elliptic curve $C$ and a certain Diophantine condition is violated  (see ~\ci[sect. 5.4]{GS3} for  construction of counter-examples). Arnold's theorem was extended by Ilyashenko  and Pjartli~\ci{MR549623} to the case when $C$ is the product of finitely many elliptic curves and $N_C$ is the direct sum of line bundles, and by the authors~\ci{GS4} for a complex torus with a Hermitian flat $N_C$ satisfying a Diophantine condition. Of course, the study of neighborhood of embedded compact complex manifolds has a long history. The reader is referred to~\ci{GS3} and references therein on neighborhoods of compact manifolds. See also recent work of Hwang~\cite{hwang-annals}, Koike~\ci{koike-fourier}, and Loray-Thom-Touzet~\ci{loray-moscou}.

The paper is organized as follows.

In Section 2, we construct a formal normal form for holomorphic neighborhoods of $C$ that is realized as a subset of
$
 \cL H^1(T_CM)$.
 Our injectivity assertion in \rt{m-finite} remains true if  $N_C$ is merely $1$-positive
 and $H^0(C,T_CM\otimes S^\ell N_C^*)=0$ for all $\ell>1$. 
In Section 3, we find a formal normal form for \emph{tangential}  foliations of neighborhoods of $C$ that contain $C$ as a  leaf.
In Section 4, we find a formal normal form for transverse   foliations of neighborhoods of $C$.   In Section 5, we apply the theorems  in~\cite{MR0137127, MR627752} to show \rt{m-finite} and analogous classification for transverse foliations. In Section 6, we use a theorem of Camacho-Movasati-Sad~\ci{MR1967036} to show that when the genus of the compact  Riemann surface $C$ is bigger than one, there are neighborhoods of $C$ that are not linearizable. Therefore, the equivalence classes contain at least two elements.

{\bf Acknowledgments.}  Part of work was carried out when X.~G. was supported by  CNRS and UCA for a visiting position at UCA.

\setcounter{thm}{0}\setcounter{equation}{0}
	\section{A formal normal form to classify neighborhoods}


 We have mentioned Grauert's formal principle asserting two neighborhoods are holomorphic equivalent if they are formally equivalent.
To further motivate our results, let us describe  the following problems about the classifications of neighborhoods.  We recall the \emph{Kuranishi problem} mentioned by Morrow-Rossi~\ci{MR627765} for the study of neighborhoods of $C$, which is to
 construct a {\it moduli space} or a parametrization that classify neighborhoods of $C$ completely  and understand the structures of the moduli space.
  This paper provides partial answers to this problem. We will also study the analogous problems for the foliations of neighborhoods of $C$.

In this section, we will construct a formal normal form for  formal equivalence classes of holomorphic neighborhoods. Certain features of the normal form will be described in details as they will be useful in the convergent proof.

To study classifications, we will use   transition functions 
for various vector bundles in coordinate charts. Let $C$ be a compact complex manifold embedded in a complex manifold $M$.    %
We cover a neighborhood of $C$ in $M$ by open sets $V_j$ and choose  coordinate charts
 $(z_j,w_j)$ on $V_j$    %
  for $M$
 such that
$$
U_j:=C\cap V_j=\{w_j=0\}.
$$
 To have Leray coverings, we may assume for instance that $U_j,V_j$ are biholomorphic to  polydiscs.
Let ${\mathcal U}=\{U_i\}$ be a finite open covering of $C$ with coordinate
charts $z_i=\var_i(p)=(z_i^1,\dots, z_i^n)$
 defined on $U_i$.
Let
\beq\label{transitionC}
z_k=
\var_{kj}
(z_j)=
\var_k\var_j^{-1}(z_j)
\nonumber
\eeq
be the transition function of $C$ on $U_{kj}:=U_k\cap U_j$. Thus the neighborhood $M$   has transition functions on  $V_{kj}:=V_j\cap V_k$
\eq{transitionN-}\nonumber
\Phi_{kj}\colon
\begin{array}{rcl}
z_k &= & \Phi^h_{kj}(z_j,w_j):=
\var_{kj}(z_j)+l_{kj}(z_j)w_j+\phi^h_{kj}(z_j,w_j),\vspace{.75ex}
\\
w_k &= &\Phi^v_{kj}(z_j,w_j):=
 t_{kj}(z_j)w_j+\phi^v_{kj}(z_j,w_j).
\end{array}
\eeq
Here, $\phi^h_{kj}$ (resp. $\phi^v_{kj}$) are holomorphic functions that vanish to order $\geq 2$ along $w_j=0$:
\eq{vord2}\nonumber
\phi_{kj}^h(z_j,w_j)=O(|w_j|^2), \quad \phi_{kj}^v(z_j,w_j)=O(|w_j|^2).
\eeq
For abbreviation, we  also  call $\Phi:=\{\Phi_{kj}\}$ a neighborhood of $C$
and $O(2)$ denotes such a $\phi_{kj}$.

 Note that the transition functions of $N_C$ are
\eq{Nkjh}
N_{kj}(h_j,v_j)=(\var_{kj}(h_j),t_{kj}(h_j)v_j)
\eeq
and the transition functions of $T_{C}M:=(TM)|_C$ are
\beq\label{transitionTCM}
g_{kj}:=\left(\begin{matrix}s_{kj}& \ell_{kj}\\0&
t_{kj}\end{matrix}\right)
(h_j)\quad \text{on } U_j\cap U_k
\eeq
for some $n\times d$ matrix $l_{kj}$,
while $s_{kj}$ is the Jacobian matrix of $\var_{kj}$. When $l_{kj}=0$,   $T_C M$ splits as $TC\oplus N_C$ (see
 \cite[Prop.~2.9]{MR627765}). With $g_{kj}$,  $T_CM$ has a basis $\{ \tilde e_k^\mu;\mu=1,\dots,n+ d\}$. Here
\eq{ekhv}
e_k^h=\{\tilde e_k^\tau\colon \tau=1,\dots,n\},
\quad e_k^v=\{\tilde e_k^{n+\nu}\colon\nu=1,\dots, d\}
\eeq
with $e_k^v$ being the basis of $N_C$ on $U_k$.  To study the foliations, we will choose $e_k^v$ to be a flat basis when $N_C$ is \emph{flat}, i.e. its transition functions are locally constant.

Throughout the paper, we fix
\eq{Nkjad}\nonumber
 N_{kj}^*(z_j,w_j)=(\var_{kj}(z_j)+\ell_{kj}(z_j)w_j,t_{kj}(z_j)w_j).
\eeq
We emphasize that $\{N^*_{kj}\}$ does not necessarily define a neighborhood of $C$ as it may not be a cocycle. Thus, we introduce the following.
\begin{defn}\label{admN} Fix $T_CM$.
Fix a holomorphic neighborhood $N^1=\{N^1_{kj}\}$ of $C$ such that $N^1_{kj}=N^*_{kj}+O(2)$.
 If $T_CM$ splits, we always take   $N_{kj}^*=N_{kj}$   defined by \re{Nkjh}.
Let $\cL M$ be the set of  germs of holomorphic neighborhoods $\Phi$ of $C$  in complex
manifolds $\tilde M$ such that $T_CM=T_C\tilde M$ and $\Phi=N^1+O(2)$.
  \end{defn}
Throughout the paper, we assume that $\cL M$ is non-empty.  Therefore the existence of   $N^1$ is ensured.

For the holomorphic  (resp. formal) equivalence of two neighborhoods of $C$, we restrict to biholomorphic   (resp. formal) mappings   $F$  from
$\Phi$ to $\tilde\Phi$ fixing $C$ pointwise. Then $F$ can be expressed as
\begin{equation}\label{change}
F_k: \begin{array}{rcl}
z_k &= & F^h_{k}(h_k,v_k):=h_k+f^h_{k}(h_k,v_k),\vspace{.75ex}
\\
w_k &= & F^v_{k}(h_k,v_k):=v_k+f^v_{k}(h_k,v_k)
\end{array}
\end{equation}
with $f_k(h_k,0)=0$ such that $F_k\Phi_{kj}
F_j^{-1}=\tilde \Phi_{kj}$ or $F\Phi F^{-1}=\tilde\Phi$ for short.   By formal equivalence, we mean that $f_{k}(h_k,v_k)$ is a formal power series in the variable $v_k$ with holomorphic coefficients in the variables $h_k\in U_k$, for some finite cover of $C$ by open sets $\{U_k\}$.
\begin{defn}
We say that $F=\{F_j\}$ defined by \rea{change} is \emph{tangent to the identity} along $C$ and write $F=Id+f$ with $f=O(2)$,  where the latter means that $f^h_{k}(h_k,v_k)$ and $f^v_{k}(h_k,v_k)$ are holomorphic functions vanishing to order $\geq 2$ at $v_k=0$.
Let ${\cL M}/{\sim}$ (resp. ${\cL M}/{\overset{f}{\sim}}$) the equivalence of \emph{holomorphic} neighborhoods of $C$ under all such biholomorphisms (resp. formal) $F$.
\end{defn}

Next, we identity sections in symmetric powers of $N_C^*$  with coordinate changes. We will associate
\eq{fkmh}
[ f_k]^m(h_k,v_k):=\left\{\left(\sum_{|Q|=m}f^h_{k, Q}(h_k)v_k^Q, \sum_{|Q|=m}f^v_{k, Q}(h_k)v_k^Q\right)\right\}\eeq
 with a $0$-th cochain   $[\tilde f]^m\in C^0(\mathcal U, T_CM\otimes S^mN_C^*)$ defined by
\eq{tfkm}
[\tilde f_k]^m(p):=\sum_{\mu=1}^{n+d}\sum_{|Q|=m}f_{k,Q}^\mu(h_k(p),v_k(p))  \tilde e_k^\mu(p)\otimes (w_k^*(p))^Q.
\eeq

  Here and in what follow, $Q\in \mathbb{N}^d$. We will also associate
 \eq{def-tilde-f}
  [ \phi_{kj} ]^m(z_{j},w_{j}):=\sum_{|Q|=m}  \phi_{kj;Q} (z_{j}) w_{j}^Q
 \eeq
  with the $1$-cochain $\{[\tilde\phi_{kj}]^m\}\in C^1(\{U_j\},\cL O(T_CM\otimes S^m(N_C^*)))$
 defined by
 \ga\label{def-tilde-f+}
 [\tilde\phi_{kj }]^m(p)= \sum_{\mu=1}^{d+d}\sum_{|Q|=m} \phi_{{kj}; Q}^\mu (h_k(p),v_k(p))\tilde e^\mu_{k }(p)\otimes (w_{j}^*(p))^Q.
 \end{gather}
Both associations are   $\C$-linear,  one-to-one and onto. By abuse of notation, we  drop tildes in $\tilde\phi_{kj},\tilde f_k$ and we interchange \re{fkmh} and \re{tfkm} (resp.  \re{def-tilde-f} and \re{def-tilde-f+}) for computation as we wish.


We also need to identity the coboundary operator $\del$ for cochain \re{def-tilde-f+} with a coboundary operator for \re{def-tilde-f}.
By Lemma 2.7 in~\ci{GS3}, applied to $E=T_CM$, we  have 
$\widetilde{\del[f]^m}=[\tilde\phi]^m$ being equivalent to $\del[f]^m=[\phi]^m$. 
Writing in column vector $[f_{ij}]^m:=([f_{ij}^1]^m, \dots, [f_{ij}^{n+d}]^m)^t$  and recalling \re{transitionTCM}, this reads
\al\label{dfphi}
(\del [f]^m)_{ij}(z_j,v_j)&:=
g_{ij}[f_j]^m-[f_i^m]\circ N_{ij}
=[\phi_{ij}]^m(z_j,v_j). 
\end{align}

%

To construct a formal normal form for ${\cL M}$, recall that   $\dim H^q(C,V)<\infty$ for $q>0$ and any holomorphic vector bundle $V$ on   $C$;  see \ci[Thm.~3.20 and Cor., p.~161]{MR2109686}.
\begin{defn}\label{choose-base}Let $C$ be a compact complex manifold and fix a holomorphic vector bundle $N_C$ on $C$. Fix a basis
$[e^m]=([e_1^m], \dots, [e_{k_ m }^m])$ for
 $H^1(C, T_CM\otimes S^ m  N_C^*)$ that is not zero. We also fix a representative $e_i^m\in Z^1(\cL U, T_CM\otimes S^ m  N_C^*)$ for $[e_i^m]$.
Define $c^m\cdot e^m:=\sum c^m_ie^m_i$ for $c^m_i\in\C$ and $(c^2,\dots c^m)\cdot(e^2,\dots, e^m)=\sum_{j=2}^m c^j\cdot e^j$. For convenience, if $H^1(C, T_CM\otimes S^ m  N_C^*)=0$,  set $c^m=0$ and $c^m\cdot e^m=0$.
\end{defn}

\le{the-lem-} Fix $m>1$. Let $\Phi,\tilde \Phi\in \cL M$ and let $F=I+f$ be a formal mapping satisfying $f=O(m)$.
\bpp
\item If $\tilde\Phi-\Phi=O(m)$, then $[\tilde\Phi]^m-[\Phi]^m\in Z^1(\mathcal U,T_CM\otimes S^m N_C^*)$; in particular,
\eq{identify}\nonumber
[\Phi-N^1]^2\in Z^1(\mathcal U,T_CM\otimes S^2 N_C^*),\quad \forall \Phi\in\cL M.
\eeq
\item If  $
F\Phi F^{-1}=\tilde\Phi+O(m)$ holds if and only if
$ \delta [f]^m=
[\Phi]^m
-[\tilde\Phi]^m.
$
\item If   $\tilde \Phi=\Phi+O(m)$, there exist a unique $c^m\cdot e^m$ and some $\tilde F=I+O(m)$ such that $\tilde F\tilde\Phi \tilde F^{-1}=\Phi+c^m\cdot e^m+O(m+1)$.
\epp

\ele
\begin{proof}$(a)$ Recall that $N_{kj}^*(h_j,v_j)=(\var_{kj}(h_j)+\ell_{kj}(h_j)v_j,t_{kj}(h_j)v_j)$ and
\eq{oldphi}\nonumber
\Phi_{kj}(h_j,v_j)=N^*_{kj}(h_j,v_j)+\sum_{\ell=2}^m[\phi_{kj}]^\ell(h_j,v_j)+O(m+1).
\eeq
The Jacobian of $N^*_{kj}$ at $(h_j,v_j)$ applied to $(\tilde h_j,\tilde v_j)$ is given by
$$
DN_{kj}^*(h_j,v_j)(\tilde h_j,\tilde v_j)=(s_{kj}(h_j)\tilde h_j+\ell_{kj}(h_j)\tilde v_j
+\partial_{h_j}\ell_{kj}(z_j)\tilde h_jv_j,t_{kj}(h_j)\tilde v_j+\partial_{h_j}t_{kj}(h_j)\tilde h_jv_j).
$$
If functions $(\tilde h_j,\tilde v_j)=O(|v_i|)$, we simplify it as
$$
DN_{kj}^*(h_j,v_j)(\tilde h_j,\tilde v_j)=DN_{kj}^*(h_j,0)(\tilde h_j,\tilde v_j)+O(|v_i|^2).
$$
Let us use transition functions $g_{kj}$ of $T_CM$ to write
$$
DN_{kj}^*(h_j,0)(\tilde h_j,\tilde v_j)=g_{kj}(h_j)(\tilde h_j,\tilde v_j)=(s_{kj}(h_j)\tilde h_j+\ell_{kj}(h_j)\tilde v_j,t_{kj}(h_j)\tilde v_j).
$$
 Then
\begin{align}\nonumber
\Phi_{kj}\Phi_{ji}(h_i,v_i) - N_{kj}^*N_{ji}^*(h_i,v_i)&=[\phi_{kj}]^{\leq m}(\var_{ji}(h_i)+\ell_{ji}(h_i)v_i,t_{ji}(h_i)v_i)
\\
&\quad+DN_{kj}^*(N_{ji}^*(h_i,v_i))[\phi_{ji}]^{\leq m}+O(m+1)\nonumber \\
&=[\phi_{kj}]^m (N_{ji}(h_i,v_i))
+g_{kj}(\var_{ji}(h_i))[\phi_{ji}]^m(h_i,v_i)\nonumber \\
&\quad + R_{kji}(h_i,v_i,[\phi_{\bullet}]^2,\dots, [\phi_{\bullet}]^{m-1})+O(m+1),\nonumber\\
&=\Phi_{ki}(h_i,v_i)-N_{kj}^*N_{ji}^*(h_i,v_i)
\nonumber
\end{align}
where function $R_{kji}$ 
is independent of $[\phi_{\bullet}]^\ell$ for $\ell\geq m$.  Applying the same computation to $\tilde \Phi=\Phi+O(m)$ and subtracting them, we obtain, writing $\psi_{kj}:=\phi_{kj}-\tilde\phi_{kj}$,
$$
[\psi_{kj}]^m (N_{ji}(h_i,v_i))
+g_{kj}(\var_{ji}(h_i))[\psi_{ji}]^m(h_i,v_i)-[\psi_{ki}]^m(h_i,v_i)=0.
$$
According to \cite[Lemma 2.7]{GS3}, this is equivalent to saying that
$$
\{[\phi_{kj}-\tilde\phi_{kj}]^m\}_{kj}\in Z^1(\mathcal U,T_CM\otimes S^m N_C^*).
$$
%
%

$(b)$ Note that when $f=O(m)$ and $\tilde\Phi=F\Phi F^{-1}$, we have $\tilde\Phi=\Phi+O(m)$.  Then $[\tilde\Phi]^m=[F\Phi F^{-1}]^m$ is equivalent to $\del [f]^m=[\phi]^m-[\tilde \phi]^m$   according to \cite[Lemma 2.16 (2.34)]{GS3}.

$(c)$ By (a), we have $\del( [\tilde\Phi]^m-[\Phi]^m)=0$. By the definition of the basis $e^\ell$, we can find $[f]^m$ such that $[\Phi]^m+ c^m\cdot e^m-[\tilde\Phi]^m=\del[f]^m$ for some $[f]^m=\{[f_k]^m\}$. Set $F_k=I+[f_k]^m$. We get $F^{-1}\tilde\Phi F=\Phi+c^m\cdot e^m+O(m+1)$ by (b).
\end{proof}

To construct normal forms, we must refine the above
  degree-by-degree normalization.
Define the approximate automorphism groups
$$
{\aut}_ m  (\Phi)=\{F=I+O(2)\colon F^{-1} \Phi F= \Phi+O( m +1)\}, \quad m=2,3,\dots.
$$
Here it is important that we allow $F$ admits lower order terms.
Since $F(\Phi+O( m +1))=F(\Phi)+O( m +1)$, then $F\in\aut_ m (\Phi)$  if and only if
$$
F\Phi-\Phi F=O({ m +1}).
$$
The latter implies that $
\Phi F^{-1}=F^{-1}\Phi +O({ m +1}).
$
Thus,   ${\aut}_ m (\Phi)$ is indeed a group. The group structure implies that the conjugacy by elements in $\aut_ m (\Phi)$ induces an equivalence relation on
$$
\cL M_m(\Phi):=\{\tilde\Phi\in \mathcal M\colon \tilde\Phi=\Phi+c^ m\cdot  e^ m +O( m +1) \}.
$$
Namely, if $\Psi,\tilde\Psi\in \cL M_m(\Phi)$,  define the equivalence relation $\sim_{\aut_m(\Phi)}$ such that
$$
\Psi\sim_{\aut_m(\Phi)}\tilde\Psi
$$
if and only if there is an $F\in \aut_m(\Phi)$ such that $F\Phi F^{-1}=\tilde\Phi+O(m+1)$.
Let $\cL M_m(\Phi)/\aut_m(\Phi)$ be the set of equivalent classes. For each equivalence class, \emph{fix} a representative.
Define $\hat{\cL C}_m(\Phi)$ be the set of elements $c^me^m$ such that $\Phi(c^m):=\Phi+c^m\cdot e^m+O(m+1)$ are among the (chosen) representatives (It will be clear from the context that $\Phi(c^m)$ do not stand for the evaluation of $\Phi$ at $c^m$).
Thus, we can express the set $\cL N_m(\Phi)$ of representatives as

$$
\{\Phi(c^m):=\Phi+c^m\cdot e^m+O(m+1)\colon c^m\cdot e^m\in \hat{\cL C}_m(\Phi)\}.
$$
 To ensure stability, $\Phi$ is always a representative, i.e.
\eq{}\nonumber
0\cdot e^m\in \hat{\cL C}_m(\Phi), \quad \Phi(0)=\Phi.
\eeq

Using this equivalence relation, we  define convergent \emph{partial normal forms} $N^m$ and formal normal forms $N^\infty$ as follows.
\begin{defn}\label{the-def} Fix  transition functions $N^1=\{N^1_{kj}\}$ for all neighborhoods
 of $C$ with a given $T_CM$ as in \rd{admN}.
\bppp
\item Define $\cL N_1=\{N^1\}$, which has one element.
Let $\cL N_2(N^1)$ to be the set of representatives $N^2(c^2)=N^{1}+c^2\cdot e^2+O(3)$ for
elements in $\cL M_2(N^1)/\aut_2(N^1)$, which are determined by
\eq{}\nonumber
c^2\cdot e^2\in \hat {\cL C}_2(N^1).
\eeq
Inductively, let $\cL N_m$  be the set of representatives
 \eq{}\nonumber
N^m(c^2,\dots, e^m):=N^{m-1}(c^2,\dots, c^{m-1})(c^m)\eeq
for elements in $\cL M_m(N^{m-1}(c^2,\dots, c^{m-1}))/\aut_m(N^{m-1}(c^2,\dots, c^{m-1}))$,  which are determined by
\eq{}\nonumber
c^m\cdot e^m\in \hat C_m(N^{m-1}(c^2,\dots, c^{m-1})).
\eeq
\item Define $N^\infty(c^2,c^3,\dots)$ to be the formal transition functions such that
$$
N^\infty(c^2,c^3,\dots)=N^m(c^2,c^3,\dots, c^m)+O(m+1), \quad \forall m.
$$
Let $\cL N_\infty$ be the set of all such formal transition functions.
\eppp
\end{defn}
By definition, $N^{m-1}(c^2,\dots, c^{m-1})(0)=N^{m-1}(c^2,\dots, c^{m-1})$. Thus, we have
\begin{prop}
Assume that $H^1(C,\tcm{m})=0$ for all $m>m_0$ and $m_0$ is finite. Then all $N^\infty=N^{m_0}$ define convergent neighborhoods and $\cL N_\infty$ is finite-dimensional.
\end{prop}

Recall that
 \eq{defH=}\nonumber
 \cL H^q(T_CM):=\bigoplus_{\ell\geq2} H^q(C,T_CM\otimes S^\ell N_C^*).
 \eeq
We now prove a formal version of  \rt{m-finite}.
\begin{thm}
\label{injection}  There exists a mapping $\frak C^f$ from $\cL M/{\, \overset{f}{\sim}}$ into $\mathcal H^1(T_CM)$. Furthermore, the mapping $\frak C^f$ is injective and  there are formal mappings transforming $\Phi\in\cL M$ into $N^\infty(\Phi)\in\cL N_\infty$, provided
\ga\label{for-F1}
 \dim \cL H^1(T_CM)<\infty, \quad \text{or} \\
 \label{for-F0}\cL H^0(T_CM)=0.
 \end{gather}
\end{thm}
\begin{proof}Recall that we fix a representant
 basis $e^\ell$ for (a basis of)
   $H^1(C,T_CM\otimes S^\ell N_C^*)$. Let $\Phi=N^1+\phi$ as in
    \rd{admN}
 define a neighborhood. We have $\del[\phi]^2=0$. Applying \rl{the-lem-} there is a unique constant vector $c_0^2(\Phi)$  such that for some $ f^2=:\{f_j^2\}\in C^{  0}(\cL U,\tcm{2})$,
$$
[\Phi
-N^1]^2=c_0^2(\Phi)\cdot e^2+\del f^2.
$$
Set $F^2_j=I+f_j^2$ and $F_2=\{F^2_j\}$. Then   $F_2\Phi F_2^{-1}=N^1+c_0^2(\Phi)\cdot e^2+O(3)$.  Recall the definition
\eq{}\nonumber
\Mt_2(N^1):=\{\widetilde\Phi\in\Mt\colon\widetilde\Phi_{kj}=N^1_{kj}+c^2\cdot e^2+O(3)\}.
\eeq
In other words, we have achieved $F_2\Phi F_2^{-1}\in \Mt_2(N^1)$. We now apply the refinement. We find a unique element  $c^2\cdot e^2\in \hat{\cL C}_2(N^1)$ and $\tilde F_2\in \aut_2(N^1)$ such that
$$
\Phi_2:=\tilde F_2F_2\Phi F_2^{-1}\tilde F_2^{-1}=N^2(c^2)+O(3).
$$
We denote this $c^2$ by $c^2(\Phi)$ and we will show that $c^2(\Phi)$ is a formal invariant.

We repeat the above two-step normalization.

Inductively, using \rl{the-lem-} we find $F_\ell$ for $\ell>2$ such that
$$F_\ell\Phi_{\ell-1} F_\ell^{-1}\in\cL M_\ell(N^{\ell-1}(c^2,\dots, c^{\ell-1})).
$$
We then find $\tilde F_\ell\in \aut_\ell(N^{\ell-1}(c^2,\dots, c^{\ell-1}))$ such that
$$
\Phi_\ell:=\tilde F_\ell F_\ell\Phi_{\ell-1}F_\ell^{-1}\tilde F_\ell^{-1}=N^\ell(c^2,\dots, c^\ell), \quad  c^\ell\in \hat{\cL C}_\ell(N^{\ell-1}(c^2,\dots, c^{\ell-1})).
$$
Again, we denote $c^\ell$ by $c^\ell(\Phi)$. Let $\ell\to\infty$. Since
$$N^{\ell+1}(c^2,\dots, c^{\ell+1})=N^{\ell}(c^2,\dots, c^{\ell})+O(\ell+1)$$
it is clear that the sequence defines  formal transition functions $N^\infty(c^2(\Phi), c^2(\Phi),\dots)$,  denoted by $N_\infty(\Phi)$.

Next we want to show that $c^2(\Phi),\dots, c^\ell(\Phi)$ are uniquely determined by the equivalence class of $\Phi$ under formal changes of coordinates that are tangent to the identity.
Suppose that $\tilde\Phi=G\Phi G^{-1}$ with $G=I+O(2)$. We have $c^2(\Phi)=c^2(\tilde\Phi)$ immediately. Suppose that $c^\ell(\Phi)=c^\ell(\tilde\Phi)=c^\ell$ for $\ell<m$. Then we can find $F,\tilde F$ so that $F\Phi F^{-1}=N^{m-1}(c')+c^me^m+O(m+1)$ and $\tilde F\tilde\Phi \tilde F^{-1}=N^{m-1}(c')+\tilde c^me^m+O(m+1)$ with $c'=(c^2,\dots, c^{m-1})$. Then $K:=\tilde FG F^{-1}$ satisfies
$$N^{m-1}(c')+\tilde c^me^m+O(m+1)=K(N^{m-1}(c')+c^me^m+O(m+1))K^{-1}\in\Mt.
$$
Thus $K\in\aut_{m}(N^{m-1}(c'))$.  Now, $\tilde c^m,c^m$ are in $\cL C_{m}(N^{m-1}(c'))$ and $N^{m-1}(c')+\tilde c^me^m+O(m+1), N^{m-1}(c')+\tilde c^me^m+O(m+1)$ are equivalent module $O(m+1)$ by $K\in\aut_{m}(N^{m-1}(c'))$. The equivalence class under the conjugacy of $\aut_{m}(N^{m-1}(c'))$ can only represented by a unique element in $\hat{\cL C}_{m}(N^{m-1}(c'))$. We obtain $c^m(\Phi)=\tilde c^m(\tilde\Phi)$. This shows that $N_\infty(\Phi)=N^\infty(c^2(\Phi),\dots,)$ is well-defined on $\cL M/{\, \overset{f}{\sim}}$. Define $N_\infty(\Phi)$ to be the image under $\frak C^f$ for the equivalence class $\Phi \mod{\!\!\overset{f}{\sim}}$  of $\Phi$ in $\cL M/{\, \overset{f}{\sim}}$.

We know that $F_\ell=I+O(\ell)$. Since a subsequence of $\tilde F_\ell$ may contain terms of a fixed order, we need to verify that the sequence $\hat F_\ell:=\tilde F_\ell F_\ell\cdots \tilde F_2F_2$ still defines a formal change of coordinates as $\ell\to\infty$.   When $ H^1(T_CM\otimes S^\ell N_C^*)=0$ for all $\ell>\ell_0$, we have $N^\infty(c^2,\dots)=N^{\ell_0}(c^2,\dots, c^{\ell_0})$ and $\tilde F_\ell=I$ and $F_\ell=I+O(\ell)$ for $\ell>\ell_0$. When $\cL H^0(T_CM)=0$, we can show inductively that $\tilde F_\ell=I$ for all $\ell\geq2$. This shows that the sequence of coordinate changes define a formal transformation $F$ and $F\Phi F^{-1}=N^\infty(\Phi)$.

Finally, we show that if $c^\ell(\Phi)=c^\ell(\tilde\Phi)$ for all $\ell$, then $\tilde\Phi$ and $\Phi$ are equivalent. Indeed, we have  $F\Phi F^{-1}=N^\infty(c^2(\Phi),\dots)=\tilde F\tilde\Phi\tilde F^{-1}$. This shows that $\tilde\Phi=\tilde F^{-1}F\Phi F^{-1}\tilde F$.
\end{proof}

 \begin{rem}We remark that there exist second-order obstructions to realize an element in $H^1(C,T_CM\otimes S^2N_C^*)$ by a neighborhood of $C$. See Griffiths~\cite{MR206980} and Morrow-Rossi~\ci[Prop. 3.4]{MR627765}.
\end{rem}

Our classification is achieved under a  group of biholomorphisms  $F$ that is smaller that the whole group of biholomorphisms, by restricting to $F$ to be tangent to the identity along $C$. Therefore, the equivalence classes under non-restricted biholomorphisms might be a smaller set; in fact by a simple dilation we can further reduce the set of equivalence classes to be a \emph{compact} but possibly a non-Hausdorff set as the case of $\dim M=2$ in~\ci[p. 323]{MR627765}.
For instance, let us consider the case of $N_C$ is a line bundle. Let $t_{kj}$ be the transition functions of $N_C$. An isomorphism of $N_C$ is given by $g_k^{-1}t_{kj}g_j=t_{kj}$. Thus $g_j$ define a global holomorphic function on $C$ without zero. Since $C$ is compact,  then the function must constant. Now it is easy to see that the transition functions for the neighborhood $\Phi_{kj}(h_j,v_j)$ is transformed into $(\Phi_{kj}^h(h_j,cv_j), c^{-1}\Phi_{kj}^v(cv_j))$.

We will give examples in Section 6 for \re{for-F1} and \re{for-F0}.

\setcounter{thm}{0}\setcounter{equation}{0}
	\section{A formal normal form for tangential foliations}
In this section, we study    \emph{tangential} foliations of neighborhoods of $C$.
  A neighborhood admits a tangential foliation, denoted by $M_\tau$,  if there are $d$ holomorphic functions   $\tilde v_k^1, \dots,  \tilde v^d_k$ such that $d\tilde v^1_k\wedge\cdots\wedge d\tilde v^d_k\neq0$ while $\tilde v_k=c_k$ and $\tilde v_j=c_j$ define the same foliation on $V_k\cap V_j$ and $v_k=0$ on $C$. The set of tangential foliations of $C$ will be denoted by $\cL M_\tau$. A biholomorphism $F$ sends $M_\tau$ into $\tilde M_\tau$, if it sends leaves of $M_\tau$ into leaves of $\tilde M_\tau$. The set of such transformations $F$ that are tangent to the identity is denoted by $\cL T(M_\tau,\tilde M_\tau)$, which depends on $M_\tau,\tilde M_\tau$.

Clearly,  the equivalences of  foliations of neighborhoods of $C$ implies the equivalence of the neighborhoods. Therefore, the classification of tangential foliations is a refinement to that of neighborhoods. This should be reflected in the construction of our normal forms for tangential foliations.

\begin{defn} Let $M_\tau$ be a tangential foliation. We call $(h_j,v_j)$ are \emph{tangential coordinates} if $v_j=cst$ define the foliation and consequently the transition functions satisfy $\Phi_{kj}^v(h_j,v_j)=\Phi_{kj}^v(v_j)$. Such  transition functions, denoted by $\Phi_\tau$, are called \emph{tangential transition functions}.
\end{defn}

\le{h-trans}
\bpp
\item  A tangential foliation  $M_\tau$ admits tangential coordinates $ (h_j,v_j)$.  That $(h_k,v_k)$ are tangential coordinates for $M_\tau$ if and only if the transverse components $\Phi^v_{kj}$ of their transition functions $\Phi_{kj}$   depend only on $v_j$.
\item Suppose  $F=I+O(2)$. Then $F$ sends a tangential foliation  $\Phi_{kj}$ into another tangential foliation   $\tilde\Phi_{kj}$ if and only if $F\Phi F^{-1}=\tilde\Phi$ and $F^v_j(h_j,v_j)=F^v_j(v_j)$ depends only on $v_j$, in which case
\eq{}\nonumber
\tilde\Phi_{kj}^v=F_k^v\Phi^v_{kj}( F_j^v)^{-1}, \quad \tilde\Phi^h_{kj}=
F_k^h\Phi_{kj}F_j^{-1}.
\eeq
\epp
\ele
\begin{proof} (a) Suppose that $M_\tau$ is defined by $\tilde v_j=cst$ such that $v_j=0$ and $d\tilde v_j^1\wedge\cdots\wedge d\tilde v_j^d\neq0$ on $C$.  Without loss of generality, we may assume that $\tilde v_j^\ell(h_j,v_j)=v_j+O(2)$.  Define $\hat v_k=\hat\Phi^v_{kj}(\hat v_j)$. Consequently, for $F_j(h_j,v_i)=(h_j,\hat v_j)$, we get $F_kF_j^{-1}(h_j,\hat v_j)=(\hat\Phi_{kj}^h(h_j,\hat v_j),
\hat\Phi_{kj}^v(\hat v_j))$.

(b) Suppose that $(h_j,v_j), (\tilde h_j,\tilde v_j)$ are  tangential coordinates and $\Phi_{kj},\tilde\Phi_{kj}$ are the corresponding transition functions for foliations $M_\tau,\tilde M_\tau$. Suppose that $F$ sends $(\cL T,C)$ into $(\widetilde{\cL T},C)$. Write $F_j=(F_j^h,F_j^v)$. Then $F^v_j(h_j,v_j)$ must be constant if $v_j$ is constant, i.e. $\tilde v_j=F_j^v(h_j,v_j)=F^v_j(v_j)$. Combining it with $\tilde v_k=\tilde\Phi^v_{kj}(\tilde v_j)$ and $v_k=\Phi_{kj}^v(v_j)$, we obtain   $\tilde \Phi_{kj}^v=F_k^v\Phi^v_{kj}(F_j^v)^{-1}$. Of course, we still have $\tilde\Phi^h_{kj}=F_k^h\Phi_{kj}F_j^{-1}$.
\end{proof}

  \begin{defn}
The set of  transformations $F=I+O(2)$ satisfying $F^v(h_j,v_j)=F_j^v(v_j)$ is denoted by $\cL T_\tau$. The set of transformations satisfying additional condition $F=I+O(m)$ is denoted by $\cL T^\tau_m$.
\end{defn}

  It is known, see for instance~\ci{GS3}, that if a tangential foliation exists, then $N_C$ is flat. We fix a flat basis $\hat e_k$ for $N_C$ over $U_k$.  Recall that when we  chose  the basis $(e_k^h,e_k^v)$ in  \re{ekhv} for $T_CM$ on $U_k\subset C$, we have decided  $e_k^v$ to be a flat basis of $N_C$ when $N_C$ is flat. Let $w_k^v$ be the dual of $e_k^v$.
     Denote by
   $
   \cL C_{cst}^\ell(\om)$
    the  sections of $T_CM\otimes S^\ell N_C^*$ of the form
     $$
     \sum_{i=1}^n\sum_{|Q|=\ell} a_{i Q}(h_k)e_{k,i}^h\otimes (w^v_{k})^Q+\sum_{\all=1}^d\sum_{|Q|=\ell} b_{\all Q}e_{k,\all}^v\otimes (w_{k}^v)^Q
     $$
  where $a_{iQ}$ are holomorphic functions and $b_{\all Q}$ are locally constant functions on $\om$ if $\om$ is a subset of $U_k$.
 Then $\del$ defined by \re{dfphi} is still a coboundary operator from $C_{cst}^q$ to $B_{cst}^{q+1}$. Denote by $\check H^q(\cL U,\cL C_{cst}^\ell)$ the \v{C}ech cohomology groups.
Define
\eq{}
\check{\mathcal H}^q_\tau(T_CM):=\bigoplus_{\ell>1}\check{ H}^q(\cL U,\cL C_{cst}^\ell).
\eeq

We now start to adapt the normal form of neighborhoods of $C$ for the tangential foliations of neighborhoods of $C$.

\begin{defn}\label{admN-3.4} Fix $T_CM$ with flat $N_C$.
Fix a tangential foliation $N_\tau^1=\{N^1_{\tau,kj}\}$ of $C$ such that $N^1_{\tau, kj}=N^*_{kj}+O(2)$.
 If $T_CM$ additionally  splits, we take $N_\tau^1=N$   defined by \re{Nkjh}. \end{defn}
Throughout the paper, we assume that $\cL M_\tau$ is non-empty. Each element in $\cL M_\tau$ is given by tangential transitions functions.   Therefore the existence of   $N_\tau^1$ is ensured regardless if $T_CM$ splits or not.

\begin{defn}
Let $\cL M_\tau$ be the set of  holomorphic tangential foliations containing $C$ as a leaf. Let $\cL M_\tau/{\,\overset{f}{ \sim}}$ the set of equivalence classes under formal tangential-foliation mappings in $\cL T_\tau$.
\end{defn}

We now state a formal classification for tangential foliations.
\begin{thm}
\label{injection-h} Let $N_C$ be flat.
Assume that $\check{ H}^1_\tau(C,\cL C_{cst}^\ell)$ are finite dimensional for all $\ell>1$.
There is a mapping $\frak C^f_\tau$ from $\cL M_\tau/{\,\overset{f}{ \sim}}$ into
$\check{\mathcal H}_\tau^1(T_CM)$ such that  if
\ga\label{for-F1-tau}\nonumber
 \dim \cL H_\tau^1(T_CM)<\infty, \quad \text{or} \\
 \label{for-F0-tau}\cL H^0_\tau(T_CM)=0,\nonumber
 \end{gather}
then $\frak C^f_\tau$ is injective and  there are formal mappings in $\cL T_\tau$ transforming $\Phi\in\cL M_\tau$ into $N_\tau^\infty(\Phi)$.
In particular, if $\cL H^1_\tau(T_CM)=0$, all tangential foliations are formally equivalent.
\end{thm}
We fix a basis $e_\tau^\ell$ for $\check{ H}^1(\cL U,\cL C_{cst}^\ell)$.
The proof is almost identical to the proof for the general case. We will only give an outline below.

\le{the-lem-TAU} Fix $m>1$. Consider transformations $F=I+f\in \cL T_\tau$ with $f=O(m)$. Suppose that $\Phi_\tau,\tilde\Phi_\tau\in\cL M_\tau$ and
$$
\tilde\Phi_\tau=\Phi_\tau+O(m),\quad
\Phi_\tau=N^1_\tau+\phi, \quad \tilde\Phi_\tau=N^1_\tau+\tilde\phi.
$$
Then we have the following.
\bpp\item  $[\tilde\Phi]^m-[\Phi]^m\in Z^1(\mathcal U,\cL C_{cst}^m)$; in particular,
$$
[\Phi-N_\tau^1]^2\in Z_{cst}^1(\mathcal U,\cL C_{cst}^2),\quad \forall \Phi\in\cL M_\tau.
$$
\item
$F\in\cL T_m^{\tau}(\Phi_\tau,\tilde\Phi_\tau)$ holds if and only if
$ \delta [f]^m=[\tilde\Phi_\tau]^m-[\Phi_\tau]^m$.
\item There exist a unique $c_\tau^m\cdot e_\tau^m\in \check{ H}^1(\cL U,\cL C_{cst}^m)$ and some $f=O(m)$   such that $ \hat\Phi:=F^{-1}\tilde\Phi_\tau  F=\Phi_\tau+c_\tau^m\cdot e_\tau^m+O(m+1)$, i.e.
    $$
    \hat\Phi=\Phi_\tau+c_\tau^m\cdot e_\tau^m+O(m+1), \quad \hat\Phi\in\cL M_\tau.
    $$
\epp
\ele

For $\Phi\in\cL M_\tau$, define
\ga{}\nonumber
\Mt_m^\tau(\Phi):=\{\tilde\Phi\in\cL M_\tau\colon \tilde\Phi=\Phi+c_\tau^2\cdot e_\tau^2+O(m+1)\},
\\
\aut_m^\tau(\Phi)= \{F\in\cL T_\tau\colon F^{-1} \Phi F= \Phi+O( m +1)\}.\nonumber
\end{gather}
Then $\aut^\tau_m(\Phi)$ is a group inducing an equivalence relation on $\Mt_m^\tau(\Phi)$. Fix an element in each equivalence class. Thus the set $\Mt_m^\tau(\Phi)/\aut^\tau_m(\Phi)$ of equivalent classes can be written as
$$
\Phi(c^m_\tau)=\Phi+c_\tau^m e^m+O(m+1)\in\cL M_\tau, \quad c^\tau e^\tau\in \cL H^m_\tau(\Phi)\subset \check H^1_{cst}(\cL U,\cL C_{cst}^m).
$$
For stability, we always choose $\Phi(0)=\Phi$.
\begin{defn}\label{the-def-TAU}Let $\Phi\in\cL M_\tau$. \bppp \item By \rl{the-lem-TAU}, we find $F_2=I+O(2)\in \cL T_\tau$ such that $F_2\Phi F_2^{-1}=N_\tau^1+\tilde c_\tau^2e_\tau^2+O(3)$, i.e.
$$
F_2\Phi F_2^{-1}\in\cL M^\tau_2(N_\tau^1).
$$
Take $\tilde F_2\in\aut^\tau_2(N_\tau^1)$ such that
\eq{}\nonumber
\Phi_2:=\tilde F_2F_2\Phi F_2^{-1}\tilde F_2^{-1}=N_\tau^1(c_\tau^2)+O(3), \quad c_\tau^2 e^2_\tau\in\cL H^2_\tau(N_\tau^{1} ).
\eeq
\item
 Let $ m >2$. Find $F_m=I+O(m)\in\cL T_\tau$  such that
 $$
 F_m\Phi_{m-1}F_m^{-1}\in \cL M^\tau_m(N^{m-1}_\tau(c^2_\tau,\dots, c^{m-1}_\tau)).
 $$
 Choose $\tilde F_m\in\aut_m^\tau((N^{m-1}_\tau(c_2^\tau,\dots, c^{m-1}_\tau))$ such that
 \ga{}\nonumber
 \tilde F_mF_m\Phi_{m-1}F_m^{-1}\tilde F_m^{-1}=N^{m-1}_{\tau}(c^2_\tau,\dots, c^{m-1}_\tau)(c^m_\tau)+O(m+1), \\ c^m_\tau e^m_\tau \in \cL H^m_\tau(N_\tau^{m-1}(c^2_\tau,\dots, c^{m-1}_\tau))\nonumber
 \end{gather}
 \item Set $N^m_\tau(c^2_\tau,\dots, c^m_\tau)=N^{m-1}_{\tau}(c^2_\tau,\dots, c^{m-1}_\tau)(c^m_\tau)$. The formal normal form of $\Phi$ is $N_\tau^\infty$ with
 \eq{}
 N_\tau^\infty=N^m_\tau(c^2_\tau,\dots, c^m_\tau)+O(m+1), \quad m=1,2,\dots.\nonumber
 \eeq
 Define $(c^2_\tau e^2_\tau,\dots)\in \cL H^1_{\tau}(T_CM)$ to be $\frak C^f_\tau(\Phi)$ for the equivalence class of $\Phi$ under $\cL T_\tau$.
\eppp\end{defn}
We leave the rest of details for the proof of \rt{injection-h} to the reader.

  We refer to \cite{GS3} for a different approach to the existence of holomorphic tangential foliation when the formal obstructions are absent in a stronger sense that the normal component of $\Phi$ is formally linearizable.
 Under a small divisor condition on cohomology groups depending only on a (flat) unitary $N_C$,  a convergence for linearizing the normal components is achieved in~\cite{GS3}.
\setcounter{thm}{0}\setcounter{equation}{0}
	\section{A formal normal form of  transverse foliations}

In this section, we will use   normal forms in Section 2 to describe neighborhoods that admit transverse foliations.

By a holomorphic \emph{transverse} foliation of $(M,C)$, we mean that on a neighborhood of $C$ in $M$ there a smooth holomorphic foliation with  all leaves intersect $C$ transversely.
First-order obstructions to transverse foliations for $(M,C)$ in higher dimensions or higher codimension were considered in~\cite{MR206980,MR627765}; however they did not settle the existence of transverse foliations when formal obstructions vanish except for the case of   $\dim M=2$   mentioned earlier.
  In this section, we will obtain an existence result on transverse foliations and the classification on them under suitable conditions on $T_CM$.

A neighborhood admits a \emph{vertical} or \emph{transverse} foliation,   if there are $n$ holomorphic functions $\tilde h_k^1, \dots, \tilde h^n_k$ such that
\eq{dhdt}\nonumber
d\tilde h^1_k\wedge\cdots\wedge d\tilde h^n_k\wedge dt_k^1\wedge\cdots\wedge dt_k^d\neq0
\eeq
 while $\tilde h_k=c_k$ and $\tilde h_j=c_j$ define the same foliation on $V_k\cap V_j$ on $C$. By an abuse of notation, we still denote the vertical foliation on a neighborhood $(M,C)$ by $(\cL V,C)$, or simply by $\Phi$.
In the previous section, we have seen that for a neighborhood to admit tangential foliation, $N_C$ must be flat. A transverse foliation does not impose conditions on $N_C$ as any $N_C$ is already foliation by its fibers. It however impose a useful condition that $T_CM$ must  split~\ci{MR627765}.  The formal obstructions for transverse foliations were obtained in~\cite{MR206980, MR627765}. In this section, we obtain a normal form for transverse foliations.

\begin{defn}Let $(\cL V,C)$ be a transverse  foliation. We say that $(h_j,v_j)$ are foliated coordinates for $(\cL V,C)$ if $h_j=cst$ define the foliation and consequently the transition functions satisfy $\Phi_{kj}^h(h_j,v_j)=\Phi_{kj}^h(h_j)$, in which case we call $\{\Phi_{kj}\}$ a transverse foliation for abbreviation.
\end{defn}
\begin{defn}
Two transverse  foliations $(\cL V,C)$ in $(M,C)$, $(\tilde{\cL  V},C)$ in $(\tilde M,C)$ are equivalent by a biholomorphic mapping $F$ if $F=I+O(2)$ sends the neighborhood $(M,C)$ into $(\tilde M,C)$ and it   sends each leave of $({\cL  V},C)$ into a leaf of $(\tilde{\cL  V},C)$.
\end{defn}
\le{H-trans}
\bpp\item Let $(\cL V,C)$ be a transverse foliation. Then the foliation   admits  foliated coordinates $(h_j,v_j)$.  That $(\tilde h_k,\tilde v_k)$ are foliated coordinates for $(\cL V,C)$ if and only if   their transition functions $\Phi_{kj}$ satisfy $\Phi_{kj}^h=\varphi_{kj}(h_j)$.
\item There is a biholomorphic mapping $F=I+O(2)$ sends a transverse  foliation $\Phi$ in foliated coordinates into another transverse  foliation $\tilde\Phi$ in the foliated coordinates if and only if $F^{-1}\Phi F=\tilde\Phi$ and $F^h_j(h_j,v_j)=
    h_j$, in which case
\eq{}\nonumber
\tilde\Phi_{kj}^h=\Phi^h_{kj}=\var_{kj}(h_j), \quad \tilde\Phi^v_{kj}=F_k^v\Phi_{kj}F_j^{-1}.
\eeq
The set of all transformations $F=I+O(2)$ with $F_j^h(h_j,v_j)=h_j$ are denoted by $\cL T_\nu$.
\epp
\ele
\begin{proof} The proof is almost verbatim to \rl{h-trans}, which is leave to the reader.
%
\end{proof}

Note that the transformations $F_j=I+f_j$ with $f_j=O(2)$ that preserves transverse foliations are rather restrictive as $f^h_j=0$. This is quite different from the study of tangential foliations in the previous section whereas transformations that preserve the tangential foliations can be higher-order perturbations in both tangential and normal components. The advantage is that the formal classification of transverse foliations is almost identical to the formal classification of the neighborhoods.

Define
$$\mathcal H_\nu^q(T_CM):=\bigoplus_{\ell>1}H^q(C,N_C\otimes S^\ell N_C^*).$$
Using \rl{H-trans}, we can obtain the following formal normal form.
\begin{thm}
\label{injection-nu}   Let $\mathcal M_\nu$ be the set of holomorphic transverse foliations of $C$.
There is a mapping $\frak C^f_\nu$ from $\cL M_\nu/{\,\overset{f}{ \sim}}$ into $\cL H_\nu^1(T_CM)$ such that if
\ga\label{for-F1-nu}\nonumber
 \dim \cL H^1_\nu(T_CM))<\infty, \quad \text{or} \\
 \label{for-F0-nu}\cL H^0_\nu(T_CM)=0,\nonumber
 \end{gather}
then $\frak C^f_\nu$ is injective and  there are formal mappings in $\cL T_\nu$ transforming $\Phi\in\cL M_\nu$ into $N^\nu_\infty(\Phi)$.
In particular, if $\cL H^1_\nu(T_CM)=0$, all transverse foliations are formally equivalent.
\end{thm}
\begin{proof}Given a transverse foliation, we know that $T_CM$ splits.
Let $\Phi\in\cL M_\nu$. Define
$$
\aut_m^\nu(\Phi)=\{F\in \cL T_\nu\colon F\Phi F^{-1}=\Phi+O(m+1)\}.
$$
 By conjugacy the group yields an equivalence relation on
$$
\cL M_m(\Phi):=\{\Phi+c^m_\nu e^m_\nu+O(m+1)\}\cap \cL M_\nu.
$$
Select representatives for the equivalence classes and denote the set of corresponding elements $c^m_\nu e^m_\nu$ by
$
\cL H^\nu_m(\Phi).
$
Thus the set of equivalence classes is given by
$$
\cL N_m(\Phi)=\{\Phi(c^m_\nu)=\Phi+c^m_\nu e^m_\nu+O(m+1)\colon c^m_\nu e^m_\nu\in \cL H_m^\nu(\Phi)\}\subset\cL M_\nu.
$$

 By \rl{H-trans},  we find $F_2=I+O(2)\in\cL T_\nu$
 such that $F_2\Phi F_2^{-1}=N^1_\nu+c_\nu^2e_\nu^2+O(3)$
 with $N^1_\nu=N^1$. Thus
$$
F_2\Phi F_2^{-1}\in\cL M^\nu_2(N_\nu^1).
$$
Take $\tilde F_2\in\aut^\nu_2(N^1)$ such that
\eq{}\nonumber
\Phi_2:=\tilde F_2F_2\Phi F_2^{-1}\tilde F_2^{-1}=N_\nu^1(c_\nu^2)+O(3), \quad c_\nu^2 e^2_\nu\in\cL H^2_\nu(N_\nu^{1} ).
\eeq
Set $N^2_\nu(c^2_\nu e^2_\nu):=N^1_\nu(c^2_\nu e^2_\nu)$.

 Let $ m >2$. Find $F_m=I+O(m)\in \cL T_\nu(\cL M_\nu)$ such that
 $$
 F_m\Phi_{m-1}F_m^{-1}\in \cL M^\nu_m(N^{m-1}_\nu(c_2^\nu,\dots, c^{m-1}_\nu)).
 $$
 Choose $\tilde F_m\in\aut_m^\nu((N^{m-1}_\nu(c_2^\nu,\dots, c^{m-1}_\nu))$ such that
 \ga{}\nonumber
 \tilde F_mF_m\Phi_{m-1}F_m^{-1}\tilde F_m^{-1}=N^{m-1}_{\nu}(c^2_\nu,\dots, c^{m-1}_\nu)(c^m_\nu)+O(m+1), \\ c^m_\nu e^m_\nu \in \cL H^m_\nu(N^{m-1}(c^2_\nu,\dots, c^{m-1}_\nu)).\nonumber
 \end{gather}
  Set $N^m_\nu(c^2_\nu,\dots, c^m_\nu)=N^{m-1}_{\nu}(c^2_\nu,\dots, c^{m-1}_\nu)(c^m_\nu)$. The formal normal form of $\Phi$ is $N_\nu^\infty$ with
 \eq{}\nonumber
 N_\nu^\infty=N^m_\nu(c^2_\nu,\dots, c^m_\nu)+O(m+1), \quad m=1,2,\dots.
 \eeq
Define $\frak C^f_\nu(\Phi)$ to be $c^2_\nu e^2_\nu,\dots)\in \cL H^1_{\nu}(T_CM)$ for equivalence class of $\Phi$ under $\cL T_\nu$.

We have defined the normal forms $N_\nu^\infty(c^2_\nu,\dots)$ for $\Phi\in \cL M_\nu$. The rest is similar to the proof of \rt{injection}. We leave the details to the reader.
\end{proof}

For the rest of the section, we deal with the existence of transverse foliations using the normal form for the neighborhoods in Section 2. Let us first improve \rd{the-def} as follows.
\begin{defn}\label{the-def+}In \rd{the-def}, we select $N^\ell(c^2,\dots, c^\ell)$ satisfying
$$(N^\ell_{kj})^h(c^2,\dots, c^\ell)=\var_{kj}(h_j),
$$ if the set of neighborhoods $\Phi$ satisfies
$$
\Phi=FN^{\ell-1}(c^2,\dots, c^{\ell-1})F^{-1}+O(\ell)
$$
for some $F\in\aut_{\ell-1}(N^{\ell-1}(c^2,\dots, c^{\ell-1}))$ has a neighborhood that admits a holomorphic transverse foliation.
\end{defn}
\begin{prop}\label{existence-formal}Let $\ell_0>1$ be an integer. Assume that
\ga{}\nonumber
H^1(C,T_CM\otimes S^\ell N_C^*)=0,\quad  \forall\ell>\ell_0.
\end{gather}
 Then $\Phi$ admits a formal transverse foliation if and only if there is a formal mapping $F\in\cL T$ such that
$F\Phi F^{-1}=N^{\ell_0}(c^2,\dots, c^{\ell_0})+O(\ell_0+1)$ with $(N^{\ell_0}(c^2,\dots, c^{\ell_0}))_{kj}^h=\varphi_{kj}(h_j)$.
\end{prop}

\begin{prop}Assume that $H^0(C, TC\otimes S^\ell N_C^*)=0$ for all $\ell>1$.   Then  two transverse foliations of neighborhoods of $C$  are equivalent if (and only if) the neighborhoods are equivalent.
\end{prop}
\begin{proof}
Suppose that $F=I+O(2)$ sends a transverse foliation $\Phi$ of a neighborhood $M$ into a transverse foliation $\tilde\Phi_{kj}$ of a neighborhood $\tilde M$. We have $\tilde\Phi_{kj}F_j=F_k\Phi_{kj}$.
Suppose that $F_j^h=I+[f^h_j]^m+O(m+1)$. We want to show that $[f^h_j]^m=0$. We know that $\Phi_{kj}^h(h_j,v_j)=\phi_{kj}(h_j)$ and $\tilde\Phi^h_{kj}(\tilde h_j,\tilde v_j)=\var_{kj}(\tilde h_j)$. We have
$$
\var_{kj}(I+[f_j^h]^m+O(m+1))=\var_{kj}(h_j)+[f_k^h(\Phi_{kj})]^m+O(m+1).
$$
Collecting terms of order $m$ in $v_j$, we obtain   $\del^h\{[f^h]^m\}=0$. Thus $[f^h]^m$ is a global section of $TC\otimes S^m N_C^*$. We conclude $[f^h]^m=0$.
\end{proof}

When $C$ is a compact Riemann surface with genus $g$, $H^0(C, TC\otimes S^\ell N_C^*)=0$
if $\deg N_C>\max\{0,g-1\}$ and $\ell>1$; see Section 6.

%

\setcounter{thm}{0}\setcounter{equation}{0}
	\section{convergence of two  classifications   and criteria for transverse foliations}
In this section, we establish the convergence results for normalizations of full neighborhoods and transverse foliations by  applying    convergent results of Grauert~\ci{MR0137127}, Griffiths~\cite{MR206980}, and Commichau-Grauert~\ci{MR627752}. We also obtain a criteria for the existence of transverse foliations using our normal forms.

  The first-order obstructions for the existence of transverse foliations were studied in~\cite{MR206980}, Morrow-Rossi~\ci{MR627765}, where the convergence of transverse foliations are not addressed with the exception that $\dim M=2$ in~\ci{MR627765} for which the classification of the foliations was   addressed.

According to Grauert~\cite[Def. 1, p.~342]{MR0137127}, we say that $E$ is weakly negative if its zero section has a relatively compact strictly pseudoconvex neighborhood. We say that $E$ is weakly positive if $E^*$ is weakly negative~\cite[Def.~2, p.~342]{MR0137127}.

Grauert~\cite[Satz 1, p.~341]{MR0137127} proved that if the zero section of a vector bundle $E$ is exceptional, then $E$ is weakly negative.
According to Grauert~\cite[Def. 3, p.~339]{MR0137127},  a connected compact complex manifold $C$ of positive dimension  in $M$ is called \emph{exceptional}, if there is a complex manifold $M'$ and a proper surjective holomorphic map $\phi\colon M\to M'$ such that $\phi(C)$ is a point $p$, $\phi\colon X\setminus C\to M'\setminus\{p\}$ is biholomorphic. In~\cite[Thm. 5,
 p.~340]{MR0137127}, Grauert proved that $C$ is  exceptional, if $C$ has a basis of strongly pseudoconvex neighborhoods.
It was proved by Grauert for codimension one compact complex manifold $C$ and by Hironaka and Rossi for higher codimension that the formal principle holds for exceptional sets.

Following Commichau-Grauert~\ci{MR627752}, we say that a vector bundle $V$ on $C$ is  \emph{$1$-positive}, if there is a  tube neighborhood $W$ of the zero section $C$ of $V$ such that the Levi-form of $V$ has at least $1$-negative eigenvalue and $W\cap V_x$ is star-shaped for each $x\in C$ and $V_x$ intersects $\pd W$ transversely. We say that $V$ is \emph{$q$-positive}, if the $\pd W$ has $q$ negative Levi eigenvalues.
We recall a lemma   mentioned by Griffiths~\cite{MR206980}.
\begin{lemma}[\!\!\cite{MR206980}, Lemma 2.1]\label{griffiths}
Let $\pi\colon N_C\to C$ be the normal bundle of $C$.    Let $W$ be a  neighborhood of $C$ in $N_C$. Let $F$ be any holomorphic vector bundle on $C$ and let $\pi^*F|_W$ be the pull-back bundle on $W$. Then
\eq{finite-dim}
\sum_{\ell\geq0}\dim H^q(C,F\otimes S^\ell N_C^*)\leq\dim H^q(W,\pi^*F|_W).
\eeq
\ele
\begin{proof} We provide a proof   using our notation for $q=1$ only.
 Let $\cL U$ be a   covering of $C$ by open sets $U_j$. We may assume that there are open sets $V_j$ in $W$ so that $U_j=V_j\cap C$. Using additional open subset of $W$, we assume that $\{V_j\}$ is a covering $\cL V$ for $W$. Furthermore, we may assume that $U_j$ and $ V_k$ are Stein,  $E$ and $ F$ holomorphically trivial on $U_j$.

The vector bundles $F$ and $N_C$ play  different roles. For $F$, we simply pull it back as $\pi^*F$ as vector bundle on $N_C$. This allows turn sections of   $F\otimes S^\ell N_C^*$ into $\pi^*F$ valued homogeneous polynomials on the manifold $N_C$.  To be specific, let us identify a section of $F\otimes S^\ell N_C^*$ by a $\C$-linear injection
$$
\iota \colon C^q(\cL U, F\otimes S^\ell N_C^*)\to C^q(\cL V, \cL O(\pi^*F|_W)).
$$
 Recall that $\pi^{-1}(U_j)$ has coordinates $(h_j,v_j)$.  Let $\{\tilde e_j^{\mu}\}$ (resp. $\{e_j^{\mu}\}$) be a basis of $F_{|U_j}$ (resp. $N_C$). Let $ \{f_{kj}\}\in C^1(\mathcal U, F\otimes S^mN_C^*)$. Then by \re{def-tilde-f+},
\eq{tfkm-M}\nonumber
 f_{kj}(p)=\sum_{\mu=1}^{n+d}\sum_{|Q|=m}f_{kj,Q}^\mu(h_k(p),v_k(p))  \tilde e_k^\mu(p)\otimes (w_j^*(p))^Q.
\eeq
Let $\tilde p\in \pi^{-1}(p)$ and $p\in U_j\cap U_k$. In coordinates, we have    $\tilde p=(p,e_j(p)\cdot   v_j(p))$. Define
$$
(\iota f)_{kj}(\tilde p)=\sum_{\mu=1}^{n+d}\sum_{|Q|=m}f_{kj,Q}^\mu(h_k(p),  v_k(p)) v_j^Q(p)  \tilde e_k^\mu(p).
$$
Now $\{\pi^{-1}(U_j)\}$ is an open covering $\hat{\cL V}$ of $F$ as a complex manifold and $\iota f\in C^{1}(\hat{\cL V}, \cL O(F))$. Let $\tilde{\cL V}$ be the covering $W$ defined by $\{V_i\cap\pi^{-1}U_j\}$. Note that $V_i\cap\pi^{-1}U_j$ are still Stein.  Analogously, if $ \{u_j\}\in C^0(\mathcal U, F\otimes S^mN_C^*)$, then
$$
 u_{k}(p)=\sum_{\mu=1}^{n+d}\sum_{|Q|=m}u_{k,Q}^\mu(h_k(p),v_k(p))  \tilde e_k^\mu(p)\otimes (w_k^*(p))^Q.
$$
 Define
$
(\iota u)_{k}(\tilde p):=\sum_{\mu=1}^{n+d}\sum_{|Q|=m}u_{k,Q}^\mu(h_k(p),v_k(p)) v_k^Q(p)  \tilde e_k^\mu(p).
$
  In the other way around, given such en expression $U_k(q)$, we associate a 0-cochain $ \tilde U=\{\tilde U_j\}\in C^0(\mathcal U, F\otimes S^mN_C^*)$.
Then $\del \iota =\iota\del$.

We want to show that $\iota$ induces an injection
$$
\iota\colon \check H^{  1}(\cL U,F\otimes S^\ell N_C^*)\to \check H^{  1}(\tilde{\cL V}, \cL O(\pi^*F|_W)).
$$
 Suppose that $f\in Z^{  1}(\cL U,F\otimes S^\ell N_C^*)$ and $\iota f=\del u$ with $u\in C^{0}(\tilde{\cL V}, \cL O(\pi^*F|_W))$. If $V_j\cap C$ is non-empty,  we have $V_j\cap C=U_j$ and we can expand
$
u=\sum u_{m},
$
where   $u_m$ is a homogeneous polynomial in $v_j$ of degree $m$.
 Then  $f=\del\tilde u_\ell$. This gives us \re{finite-dim}.
\end{proof}
It is known from the Andreotti-Grauert
theory~\ci{MR150342}  that if $\pd W$ has $(q+1)$ negative Levi eigenvalue or $(n+d-q)$ positive Levi eigenvalues at each boundary point, then
$$
\dim H^q(W,\cL F)<\infty
$$
for any coherent analytic sheaf $\cL F$ on $W$.
Our convergent classifications are based on the following theorems, which can be regarded as strong forms of the formal principle.
\begin{thm}[\!\!\cite{MR0137127}, Satz 7, p.~363]
\label{conv1} Let $N_C$ be negative. Then there exists an integer $\ell_0$ such that
$H^1(C,T_CM\otimes S^\ell N_C^*)=0$ for $\ell>\ell_0$. Furthermore, if  $\hat F_j$ are holomorphic mappings such that  $\tilde\Phi_{kj}=\hat F_k^{-1}\Phi \hat F_j+O(\ell_0+1)$, then there are holomorphic mappings $F_j$ such that $F_j=\hat F_j+O(\ell+1)$ and $F_k^{-1}\Phi_{kj}F_j=\tilde\Phi_{kj}$.
\end{thm}
\begin{rem}\label{HR-remark}
Let $C$ be a compact Riemann surface and let $N_C$ be a line bundle with $\deg N_C<0$,
By Riemann-Roch (see \re{expansion}-\re{expansion+} below),  we have
$$\lim_{\ell\to\infty}\ell^{-1}\dim H^0(C,T_CM\otimes S^\ell N_C^*)=-\deg N_C.$$ Using this one can prove that there are divergent formal mappings that preserve the germs of neighborhoods of the zero section of $N_C$; see \rp{MR-issue} below.
\end{rem}

On the other hand, with $1$-positivity, we have the following {\em strong } formal principle.
\begin{thm}[\!\!\cite{MR627752}, Satz 4, p.~119] \label{conv2}Let $N_C$ be $1$-positive. Let $ F_j$ be formal biholomorphic mappings such that $\tilde\Phi_{kj}=F_k^{-1}\Phi_{kj}F_j$. Then $F_j$ must converge.
\end{thm}
Commichau-Grauert~\ci{MR627752} proved the theorem for codimension $1$ case first.  Their proof for the  higher codimensions follows from a blow-up  ~\ci[p.~115 and p.~126]{MR627752} along $C$. Indeed,  one can verify easily as follows. By blowing up a neighborhood $M$ along $C$, we obtain $\tilde M$ and $\tilde C$ such that $\tilde C$ has codimension $1$ while $\pd\tilde M$ is biholomorphic to $\pd M$. If $F\colon (M,C)\to(M',C)$ is a formal mapping that is tangent to the identity along $C$, then the blow-up induces a formal mapping $\tilde F$ from $(\tilde M,\tilde C)$ to $(\widetilde {M'},\tilde C)$ (that may not be tangent to the identity along $\tilde C$). However, the theorem for codimension $1$ case, which is proved for any formal mapping that preserves $C$, implies that $\tilde F$ is convergent. Consequently, $F$ is also convergent.

 We now prove the convergence classification for the neighborhoods.
\begin{cor}\label{full-fd} Let $N_C$ be weakly negative or $2$-positive. Then the set of  holomorphic equivalence classes of neighborhoods of $C$ is identified with a subset of finite dimensional space $\cL H(T_CM)<\infty$. Furthermore, a neighborhood $\Phi$ admits a transverse foliation if and only if its normal formal $N^\infty(\Phi)$ satisfies
$(N_{kj}^\infty(\Phi))^h=\var_{kj}$.
\end{cor}
\begin{proof}Let us summarize the proof of \rt{injection}.
Let $\Phi$ be a holomorphic (convergent) neighborhood of $C$ such that $\Phi=N^1+O(2)$, where $N^1$ is a convergent neighborhood. Then we find a (convergent) biholomorphism $F_2$ such that $F_2^{-1}\Phi F_2=N^1+c^2(\Phi)e^2+O(3)$ and  choose a convergent neighborhood $N(c^2(\Phi))$ such that $F_2^{-1}\Phi F_2=N(c^2(\Phi))+O(3)$. Inductively, we find biholomorphism $F_2,\dots, F_m$ and convergent neighborhoods $N(c^2(\Phi)e^2, \dots, c^m(\Phi))$  such that
$$
F_m^{-1}\cdots F_2^{-1}\Phi F_2\cdots F_m=N(c^2(\Phi), \dots, c^m(\Phi))+O(m+1).
$$
Choose $m$ so that $H^1(C,T_CM\otimes S^\ell N_C^*)=0$ for $\ell>m$.
 Hence they are holomorphically equivalent by \rt{conv1}.
\end{proof}

 We obtain the classification for transverse foliations.

\begin{cor}
Let $C$ be a compact   complex manifold.
Suppose that $N_C$ is $1$ positive.
\bppp\item
  If two transverse foliations are equivalent by a formal biholomorphic mapping $F$ preserving the foliations, then $F$ is also convergent.
\item  Assume further that
$
\dim H^1(C,N_C\otimes S^\ell N_C^*)<\infty
$
for some $\ell_0$.
There is an injective mapping from the set   ${\cL M_\nu}/{\,\sim}$ of holomorphic equivalence classes of the  transverse foliations into the finite-dimensional space
$\mathcal H^1_\nu(T_CM).$
\eppp
\end{cor}
\begin{proof} Assertion $(i)$ is a consequence of \rt{conv2}. Assertion $(ii)$ follows from $(i)$ and the formal classification by \rt{injection-h}.
\end{proof}

 When $C$ is a compact Riemann surface and $N_C$ is a positive line bundle,  Morrow-Rossi, using the theorem of Commichau-Grauert, the equivalence classes of transverse foliations under foliation-preserving transformations are actually infinite dimensional.

We conclude this section by considering the classification of neighborhoods of $C$ that are not necessary tangent to the identity. We will however use coordinate changes that fix $C$ pointwise and preserve $N_C$. We consider the case of $N_C$ is a line bundle. Let $t_{kj}$ be the transition functions of $N_C$. An isomorphism of $N_C$ is given by $g_k^{-1}t_{kj}g_j=t_{kj}$. Thus $g_j$ define a global holomorphic function on $C$ without zero. Since $C$ is compact,  then the function must constant. Now it is easy to see that the transition functions for the neighborhood $\Phi_{kj}(h_j,v_j)$ is transformed into $(\Phi_{kj}^h(h_j,cv_j), c^{-1}\Phi_{kj}^v(cv_j))$. Note this non-homogenous dilation is used by Morrow-Rossi~\ci{MR627765} to get a complete set of moduli spaces when $\dim C$ and $\codim C$ are $1$, and $N_C$ is negative.
\setcounter{thm}{0}\setcounter{equation}{0}
	\section{ Riemann surfaces in complex surfaces}

We illustrate Grauert's result by showing that the obstructions exists for a neighborhood with negative normal line bundle being not holomorphically equivalent to a neighborhood of its zero section; compare  \cite[p.~211]{MR947141} on Grauert's work.
We also prove the assertion in \rrem{HR-remark}.
 This explains why we could not apply Grauert's formal principle for weakly negative $N_C$ to the classification of transverse foliations,  and  it seems  to the authors that an application of  a formal principle to the transverse foliations
 needs a statement stronger than \rt{conv1}. An interested reader should consult
Morrow-Rossi~~\ci[p.~323, line~3, and Thm. 6.3]{MR627765} for negative $N_C$.

Let us first recall some facts on the dimensions of $0$-th and first cohomology groups of line bundles.
Note that a line bundle on a compact Riemann surface is positive if and only if its degree is positive.
Let $C$ be a compact Riemann surface with genus $g$.  When $L\to C$ is a line bundle with degree $\nu_L<0$,   $L$ has no global section. Recall the duality
$$
H^q(C,\Om^p(E))=H^{1-q}(C,\Om^{1-p}(E^*))^*, \quad q=0,1.
$$
We have
$\deg TC=2-2g$ and $\deg K_C=2g-2$, where $K_C$ is the canonical line bundle   of $C$.
The Riemann-Roch theorem says that if $h^0(L)$ is the dimension of $H^0(C,\mathcal O(L))$ then
\eq{RR}\nonumber
h^0(L)-h^0(K_C\otimes L^{-1})=\deg L+1-g.
\eeq
Thus
\eq{crit}\nonumber
h^0(L)\geq \deg L-g+1.
\eeq
This provides the following useful estimates for positive $N_C$:
\ga\label{expansion}
\dim H^1(C, NC\otimes S^\ell N_C^*)\geq g-1+(\ell-1)\deg N_C,\\
\dim H^1(C, TC\otimes S^\ell N_C^*)\geq 3g-3+\ell\deg N_C.\label{expansion+}
\end{gather}

\subsection{Normal forms on neighborhoods}
\subsubsection{Negative $N_C$} In this case, the formal principle and \re{for-F1} hold.
We have
\eq{}\nonumber
H^1(C, (T_C\oplus N_C)\otimes S^\ell N_C^*)=(H^0(C, K_C\otimes(K_C\otimes S^\ell N_C+N_C^*\otimes S^\ell N_C))^*.
\eeq
Now $d_\ell^h:=\deg K_C\otimes(K_C\otimes S^\ell N_C)=4(g-1)+\ell\deg N_C$ and $ d_\ell^v:=\deg K_C\otimes(N_C^*\otimes S^\ell N_C)=2(g-1)+( \ell-1)\deg N_C$. Thus
\ga{}\nonumber\dim H^1(C, TC\otimes S^2N_C^*)\geq 3g-3+2\deg N_C,\\
\dim H^1(C, N_C\otimes S^2N_C^*)\geq d_2^v-g+1=g-1+2\deg N_C.\nonumber
\end{gather}
 By Morrow-Rossi~\ci[Thm.~6.3]{MR627765}, each element in $H^1(C,N_C\otimes S^2 N_C^*)$ can be realized by transverse foliations.
 See Camacho-Movasati-Sad~\ci{MR1967036} for a different approach.
This shows that an element corresponding to a non-zero element in $H^1(C,N_C\otimes S^2 N_C^*)$ is not equivalent to a neighborhood of the zero section of $N_C$.

\subsubsection{Positive $N_C$ and negative $T_CM\otimes S^2 N_C^*$} This occurs  if and only if $\deg N_C> \max\{0,g-1\}$, in which case \re{for-F0} holds. Then the formal principle holds. Also,   both
$$
\oplus_{\ell>1}  H^1(TC\otimes S^\ell N_C^*),\quad \oplus_{\ell>1}  H^1(C,N_C\otimes S^\ell N_C^*).
$$ are infinite dimensional.

\begin{prop}\label{realization-MR} Let $C$ be a compact Riemann surface with genus $g$. Suppose $N_C$ is a line bundle. For any finite $r$, the following hold:
\ga\label{finite-sum-case}
\oplus_{m=2}^rH^1(C,T_CM\otimes S^m N_C^*)\subset\frak C(\cL M/\sim)\subset\oplus_{m=2}^\infty H^1(C,T_CM\otimes S^m N_C^*)\end{gather}
if $ \deg N_C>\max\{0,g-1\}$; if $\deg N_C>0$ and $T_CM$ splits then
\ga{}
\oplus_{m=2}^rH^1(C,N_C\otimes S^m N_C^*)\subset\frak C_\nu(\cL M_\nu/\sim)\subset\oplus_{m=2}^\infty H^1(C,N_C\otimes S^m N_C^*).
\label{finite-sum-case-nu}
\end{gather}
\end{prop}
\begin{proof}
Our proof is based on  a  construction in \ci{MR627765} to realize $\oplus_{\ell>1}H^1(C, N_C\otimes S^\ell N_C^*)$ for transverse foliations. Thus \re{finite-sum-case-nu} is essentially in \ci{MR627765}. As indicated in \ci{MR627765}, the construction applies to
$$\oplus_{\ell>1}H^1(C, T_CM\otimes S^\ell N_C^*)$$
 as well, which we show below. We will need
 $$
 H^0(C,T_CM\otimes S^\ell N_C^*)=0, \quad\forall \ell>0
 $$
 which holds for $\deg N_C>\max\{0,g-1\}$. Thus $\aut_m(N^m(c^2,\dots, c^m))$, defined as
 $$\{F=Id+O(2)\colon FN^m(c^2,\dots, c^m)=N^m(c^2,\dots, c^m) F+O(m+1)\}$$
 consists of mappings of the form $F=Id+O(m+1)$.

   Take a holomorphic disk $U_0$ in $C$ and let $U_0'$ be a smaller disc in $U_0$. Then $U_0,U_1:=C\setminus U_0'$ form a Leray covering of $C$ as both $U_0,U_1$ are Stein. $U_0\cap U_1$ is an annulus biholomorphic to
$\{z\in\C\colon r<|z|<1/r\}$.  Since $C$ is covered by two sets, then
\eq{Z=C}\nonumber
Z^1(\cL U,L)=
C^1(U_0\cap U_1, L).
\eeq
Since a line bundle on an open Riemann surface is holomorphically trivial~\cite[p.~52]{MR0207977},  then $N_C$ is completely determined by
$$
(\var_{10}(z_0),t_{10}(z_0)v_0)
$$
where $t_{10}(z)$ is non-vanishing holomorphic function on $U_0\cap U_1$ and $\var_{10}$ is holomorphic and injective on $U_0\cap U_1$. Also $\deg N_C$ is the winding number of $t_{10}$ on the unit circle. A neighborhood of $C$ with $N_C|_{U_1}$ being trivial is precisely given by
$$
\Phi_{10}(z_0,v_0):=(\var_{10}(z_0)+l_{10}(z_0)v_0,t_{10}(z_0)v_0)+
\sum_{\ell>1}\Phi_{10,\ell}(z_0)v_0^\ell
$$
by patching
$
(U_0\sqcup U_1)\times\Del_\e/\!\sim
$
with  $(z_1,v_1)\sim\Phi_{10}(z_0,v_0)$. Here $\Phi_{10,\ell}$ are otherwise any holomorphic functions on ${U_{01}}\times \Del_\e$ subject to the condition that \eq{constr}
\sum_{\ell>1}\sup_{K}|\Phi_{10,\ell}|<C_K<\infty
\eeq
for any compact subset $K$ of $U_{01}\times\Del_\e$.  By an abuse of notation, the $(z_1,U_1)$ stands for a finite collection of   coordinate charts $(z_j,U_j)$, where $U_1=\cup \tilde U_j$ and $z_j$ maps $\tilde U_j$ onto the unit disc.

When $l_{10}=0$, $T_CM$ splits. When $\Phi_{10,\ell}^h=0$, the neighborhood admits transverse foliation. When $t_{10}, l_{10}$ and $\Phi_{10,\ell}^v$ are constant, the neighborhood admits tangential foliations.
The degree of $N_C$ is the winding number of $t_{10}$ on the unit circle.

To realize each element in \re{finite-sum-case}, we recall briefly the construction in \rt{injection}. Since $C$ is covered by two sets,  we can take $$N^1=(\var_{10}(z_0)+l_{10}(z_0)v_0,t_{10}(z_0)v_0).
$$
We now take the advantage that $\Phi_{01,\ell}$ need only to satisfy \re{constr} and are otherwise arbitrary. Fix a finite basis $e^2$ for $H^1(\cL U,T_CM\times S^2N_C^*)$. Let $e^2_j$ be represented by $C^1(\cL U, T_CM\times S^2N_C^*)$, still denoted by $e^2_j$. We also use the identification via \re{def-tilde-f} and \re{def-tilde-f+}. Thus each element $c^2e^2\in H^1(\cL U,T_CM\times S^2N_C^*)$ is associated with $c^2e^2$ in the form
$$
N^2(c^2)=(\var_{10}(z_0)+l_{10}(z_0)v_0,t_{10}(z_0)v_0)+ c^2\hat N_2^2(z_0)v_0^2.
$$
Since $C$ is covered by two open sets, $N^2(c^2)$ is indeed a cocycle for transition functions of a neighborhood of $C$.

For finitely many $c^2,\dots, c^r$, $\sum_{\ell=2}^r c^\ell e^\ell$ is associated to $N^r(c^2,\dots, c^r)$ satisfying the convergence constrain \re{constr}.
\end{proof}

\subsection{Verifying \rrema{HR-remark}}
%
%
%
%
%
%
%
%
\begin{prop}\label{MR-issue}Assume that $T_CM$ splits. Assume that $H^1(C,T_CM\otimes S^\ell N_C^*)=0$ for $\ell\geq \kappa$ for a finite $\kappa$.
\bpp\item Suppose $H^0(C,T_CM\otimes S^\ell N_C^*)\neq0$ for a sequence $\ell=\ell_j\to\infty$. Then there are divergent formal mappings that preserve the germs of neighborhoods of the zero section of $N_C$.
\item
Suppose $H^0(C,N_C\otimes S^\ell N_C^*)\neq0$ for some $\ell>1$. Then there is a possibly divergent  mapping $F$ preserving the germ of the zero section of $N_C$ while $F$ does not preserve the transverse foliation of $N_C$.
\epp
\end{prop}
\begin{proof}
 If $T_CM$ has locally constant transition functions, the proof is straightforward as each global section $[\phi]^m$ of $T_CM\otimes S^mN_C^*$ gives us an automorphism of $T_CM$ of the form $F=I+f$ where $f$ is homogeneous of degree $m$. Then $F=I+\sum_{\ell>1}f_m$ diverges for suitable choices of coefficients of the global sections of $T_CM\otimes N_
c^\ell$ for $\ell>1$.

 When $T_CM$ is not flat, we verify the assertion by a \emph{reversing} linearization procedure. Take $j$ such that $\ell_{j}\geq \kappa$.
Using a global section $f^*_{\ell_{\ell_j}}$ of $T_CM\otimes S^{\ell_j} N_C^*$, we find $F_{\ell_0}=I+f_{\ell_j}$ where $f_{k, \ell_j}$ is homogeneous of degree $\ell_j$ in $v_k$ variables. Then
\eq{m-comm}
F_{\ell_j}N^1 F_{\ell_j}^{-1}=N^1+O(\ell_j+1)
\eeq
By condition on $H^1$ and \re{m-comm}, we find a formal mapping $\tilde F_{\ell+1}=I+O(\ell_j+1)$ such that $\tilde F_{\ell_j+1}F_{\ell_j}N^1 F_{\ell_j}^{-1}\tilde F_{\ell_j+1}^{-1}=N^1$. We repeat this and find
$
F_{\ell_i}, \tilde F_{\ell_i+1}
$
for $i>j$
such that $\tilde F_{\ell_{i}+1}F_{\ell_i}$ and $N^1$ commute.
Then
$$
F_\infty=\lim_{i\to\infty} \tilde F_{\ell_i+1}F_{\ell_i}\circ\cdots\circ \tilde F_{\ell_j+1}F_{\ell_j}
$$
commutes with $N^1$. We have
$$
F_\infty=F_{\ell_i}\cdots \tilde F_{\ell_j+1}F_{\ell_j}+O(\ell_i+2),\quad
F_\infty=\tilde F_{\ell_i+1}F_{\ell_i}\cdots \tilde F_{\ell_j+1}F_{\ell_j}+O(\ell_i+2).
$$
In other words,
$$
[F_\infty]_{\ell_i}=f_{\ell_i}+R(f_{\ell_j},\dots, f_{\ell_{i-1}})
$$
When the coefficients of $\tilde f_{\ell_i}$ grows sufficiently fast as $j\to\infty$, we get a divergent $F_\infty$. This proves the first assertion.

 The second assertion can be proved by a similar argument.
\end{proof}
As an example, the condition in \rp{MR-issue} (a) is satisfied when $C$ is the Riemann sphere and $N_C$ is a negative line bundle.

\newcommand{\doi}[1]{\href{http://dx.doi.org/#1}{#1}}
\newcommand{\arxiv}[1]{\href{https://arxiv.org/pdf/#1}{arXiv:#1}}

  \def\MR#1{\relax\ifhmode\unskip\spacefactor3000 \space\fi%
  \href{http://www.ams.org/mathscinet-getitem?mr=#1}{MR#1}}

%

\begin{thebibliography}{10}

\bibitem{MR150342}
A. Andreotti and H.
   Grauert.
\newblock {Th\'{e}or\`eme de finitude pour la cohomologie des espaces complexes}.
\newblock
  {\em Bull. Soc. Math. France},
   90:
  193--259,
   1962.


\bibitem{MR0431285}
V.~I. Arnol'd.
\newblock Bifurcations of invariant manifolds of differential equations, and
  normal forms of neighborhoods of elliptic curves.
\newblock {\em Funkcional. Anal. i Prilo\v{z}en.}, 10(4):1--12, 1976.

\bibitem{MR947141}
V.~I. Arnol'd.
\newblock {\em Geometrical methods in the theory of ordinary differential
  equations}, volume 250 of {\em Grundlehren der Mathematischen Wissenschaften
  [Fundamental Principles of Mathematical Sciences]}.
\newblock Springer-Verlag, New York, second edition, 1988.
\newblock Translated from the Russian by Joseph Sz\"{u}cs [J\'{o}zsef M.
  Sz\H{u}cs].

\bibitem{MR1967036}
C.~Camacho, H.~Movasati, and P.~Sad.
\newblock Fibered neighborhoods of curves in surfaces.
\newblock {\em J. Geom. Anal.}, 13(1):55--66, 2003.

\bibitem{MR627752}
M.~Commichau and H.~Grauert.
\newblock Das formale {P}rinzip f\"{u}r kompakte komplexe
  {U}ntermannigfaltigkeiten mit {$1$}-positivem {N}ormalenb\"{u}ndel.
\newblock volume 100 of {\em Ann. of Math. Stud.}, pages 101--126. Princeton
  Univ. Press, Princeton, N.J., 1981.

\bibitem{MR3488114}
M.~Falla~Luza and F.~Loray.
\newblock On the number of fibrations transverse to a rational curve in complex
  surfaces.
\newblock {\em C. R. Math. Acad. Sci. Paris}, 354(5):470--474, 2016.

\bibitem{GS3}
X.~Gong and L.~Stolovitch.
\newblock Equivalence of neighborhoods of embedded compact complex manifolds
  and higher codimension foliations.
\newblock {\em Arnold Math J.}, (8):61–145, 2022.

\bibitem{GS4}
X.~Gong and L.~Stolovitch.
\newblock
On neighborhoods of embedded complex tori, \newblock preprint, submitted,
  \newblock https://doi.org/10.48550/arxiv.2206.06842.


\bibitem{MR0137127}
H.~Grauert.
\newblock \"{U}ber {M}odifikationen und exzeptionelle analytische {M}engen.
\newblock {\em Math. Ann.}, 146:331--368, 1962.

\bibitem{MR206980}
P.~A. Griffiths.
\newblock The extension problem in complex analysis. {II}. {E}mbeddings with
  positive normal bundle.
\newblock {\em Amer. J. Math.}, 88:366--446, 1966.


\bibitem{MR0207977}
R.~C. Gunning.
\newblock {\em Lectures on {R}iemann surfaces}.
\newblock Princeton Mathematical Notes. Princeton University Press, Princeton,
  N.J., 1966.

\bibitem{MR0171784}
H.~Hironaka and H.~Rossi.
\newblock On the equivalence of imbeddings of exceptional complex spaces.
\newblock {\em Math. Ann.}, 156:313--333, 1964.

\bibitem{MR621013}
A.~Hirschowitz.
\newblock On the convergence of formal equivalence between embeddings.
\newblock {\em Ann. of Math. (2)}, 113(3):501--514, 1981.




\bibitem{MR1152943}
J.~C. Hurtubise and N.~Kamran.
\newblock Projective connections, double fibrations, and formal neighbourhoods
  of lines.
\newblock {\em Math. Ann.}, 292(3):383--409, 1992.

\bibitem{hwang-annals}
J.-M. Hwang.
\newblock An application of {C}artan's equivalence method to {H}irschowitz's
  conjecture on the formal principle.
\newblock {\em Ann. of Math. (2)}, 189(3):979--1000, 2019.


\bibitem{MR704627}
Y.~S. Ilyashenko.
\newblock Imbeddings of positive type of elliptic curves into complex surfaces.
\newblock {\em Trudy Moskov. Mat. Obshch.}, 45:37--67, 1982.

\bibitem{MR549623}
Y.~S. Ilyashenko and A.~S. Pjartli.
\newblock Neighborhoods of zero type imbeddings of complex tori.
\newblock {\em Trudy Sem. Petrovsk.}, (5):85--95, 1979.




\bibitem{MR2109686}
 K. Kodaira.
 \newblock {\em Complex manifolds and deformation of complex structures},
   \newblock Classics in Mathematics.
   \newblock Reprint of the 1986 English edition.
   \newblock Springer-Verlag, Berlin, 2005.
   \newblock Translated from the 1981 Japanese original by Kazuo Akao.


\bibitem{koike-fourier}
T. Koike.
\newblock Linearization of transition functions of a semi-positive line bundle
  along a certain submanifold.
\newblock {\em Ann. Inst. Fourier (Grenoble)}, 71(5):2237--2271, 2021.

\bibitem{loray-moscou}
F. Loray, O. Thom, and F. Touzet.
\newblock Two-dimensional neighborhoods of elliptic curves: formal
  classification and foliations.
\newblock {\em Mosc. Math. J.}, 19(2):357--392, 2019.



\bibitem{MR1250979}
M.~B. Mishustin.
\newblock Neighborhoods of the {R}iemann sphere in complex surfaces.
\newblock {\em Funktsional. Anal. i Prilozhen.}, 27(3):29--41, 95, 1993.

\bibitem{MR627765}
J.~Morrow and H.~Rossi.
\newblock Some general results on equivalence of embeddings.
\newblock In {\em Recent developments in several complex variables ({P}roc.
  {C}onf., {P}rinceton {U}niv., {P}rinceton, {N}. {J}., 1979)}, volume 100 of
  {\em Ann. of Math. Stud.}, pages 299--325. Princeton Univ. Press, Princeton,
  N.J., 1981.

\end{thebibliography}

\end{document}